\RequirePackage{fix-cm}
\documentclass[smallextended]{svjour3}       
\smartqed  
\usepackage{graphicx,amsfonts,amssymb,amscd,amsmath,enumerate,verbatim,calc, textcomp,mathrsfs}
\usepackage{epstopdf}
\begin{document}

\title{On Coretractable Module
}


\author{Abhay Kumar Singh   \and
        Amrit Kumar Mahato 
}


\institute{Department of Applied Mathematics \at
              Indian School of Mines, Dhanbad(826004) \\
              \email{singh.ak.am@ismdhanbad.ac.in}           
           \and
           Department of Applied Mathematics \at
              Indian School of Mines, Dhanbad(826004)\\
              \email{kumaramritmaths@gmail.com} 
}

\date{Received: date / Accepted: date}

\maketitle

\begin{abstract}
The main aim of this paper is to discuss the structure of coretractble and mono-coretractble modules. We have provided the characterizations of coretractable modules in terms of Kasch ring, CS modules, mono-coretractble modules, projective modules, epi retractable modules, and retractable modules. We have also discussed some equivalent condition related with coretractbele and mono-coretractble modules.

\end{abstract}

\section{Introduction}
\label{intro}
In this paper, all rings are associative and modules are unitary right modules. Any terminology not defined hare may be found in Gooderl and Warfield \cite{1} and McConnell and Robson \cite{2}. Zelmanowitz considered compressible modules in detail in a series of papers \cite{11}.  Following \cite{12} the right R-module M is called compressible if for each non-zero zero submodule N of M, there exists a monomorphism $f: N \rightarrow M.$ Mc Connel and Robson considered the epi-retractable modules as the dual of compressible module and defined it as follows a module M is called epi-retractable if every submodule of M is a homomorphic image of M \cite{2}. An R-module M is said to be retractable if $ Hom(M, N)\neq 0 $ for any non-zero submodule N of M. Retractable modules have been discussed by some authors in series of papers 
 \cite{5}, \cite{6} ,\cite{7} ,\cite{8}, \cite{9}, \cite{10}. Dual concept to retractable modules was investigated in paper \cite{13}. B. amini et al.\cite{13} proved that the right CC rings are two sided perfect and provided condition under which free modules are coretractable. In \cite{6}, P.F. Smith introduced the concept of essentially compressible modules as a generalization of compressible modules. Follows Vedadi \cite{9}, an R-module M is said to be essentially retractable if $ Hom(M, N)\neq 0 $, for any non-zero essential submodule of M which is a natural generalization of retractable modules. In \cite{7}, Abhay K. Singh discussed some results related with essentially slightly compressible modules and rings. Some structure of completely coretractable rings have been discussed by J. Zemlicka  \cite{14}. \\
If N is a submodule of M, we write $N \leq M $and if N is an essential submodule of M then we write 
$ N \unlhd M. $ A partial endomorphism of a module M is a homomorphism from a submodule of M into M. 
 We say that N is a dense submodule of M (written $ N \subseteq  _{d}M $) if,  for any $ y \epsilon  M $ and $ x \epsilon  M / \lbrace 0 \rbrace , xy^{-1} N \neq 0 $ (i.e.,there exists $ r \epsilon R $ such that $ xr \neq  0 $ ,and $ yr \epsilon  N $ ). If $ N \subseteq _{d} M $  , we also say that M is a rational extension of N. 
 Let I = E(M),and let H= End $(I_{R}) $, operating on the left of I. We define Rational Hull of a Module $\widehat{E}(M) $ as 
$\widehat{E}(M) = \lbrace i \epsilon I : \forall h \epsilon H , h(M) = 0,  h(i) = 0\rbrace. $ 
Clearly, this is an R-submodule of I containing M.  
An R-module MR is said to be rationally complete if it has no  proper rational extensions, or equivalently $ \widehat{E}(M) = M$. \\

 A  module  is  called  mono  coretractable  if  each  of  its  factors  can  be  embedded  in  it.  That  is,  a  module  M  is  mono-coretractable if for each submodule N of M there exists a monomorphism $g: M / N \rightarrow M.$  Mono coretractable module is coretractable. Every simple and semi simple module is mono coretractable $Z_{p}^{\infty}$ as a Z-module is mono coretractable and Z as Z-module is not mono coretractable. Some charecterizations of mono-coretractable and coretractable modules discussed in \cite{15}. In section 2 we study some properties of coretractable modules. We delveloped some relation between epi-retractable modules, coretractable modules, semiprimitive modules, and Kasch ring. In section 3, we discussed the characterizations mono-retractable modules in terms of non-singular uniform modules and kasch ring.

\section{Structure of Coretractable Modules}
 A module M is called coretractable (see [13]) if $ Hom_{R}(M/K,M)\neq 0 ~\text{with}~f(K)\\=0$  for any proper submodule K of M. Module M is said to be torsion free module if $ \forall r \epsilon R ~\text{and} ~ m \epsilon M, rm = 0 \Rightarrow r=0 ~ \text{or}~ m=0. $ Before we prove some general module theoretic correction of coretractability, we investigate role of free, dense and rational submodules in coretactable        modules. In our first result we discuss a coretractable module M in terms of torsion free module.\\ \\
 Proposition 2.1. Let M be coretractable module over a commutative ring R then M is torsion free.\\ 
 Proof. Since M coretractable thus we have a homomorphism s.t.  $ Hom_{R}(M/K,M) \\ \neq 0. $ \\
 If possible let M is not torsion free then $ \forall r \epsilon R \ and \ m \epsilon M \Rightarrow rm = 0. $ \\ 
Thus we have $  f(rM/K)= f{r(m+K)} = f(rm) + f(rK) = 0, \ \text{since} \ \\ f(K) = 0.$ A contradiction to the fact that $ f\neq 0. $ Hence M is torsion free. \\ \\
 An R-module $ M_{R} $ is said to be rationally complete if it has no 
proper rational extensions, or equivalently $ \widehat{E}(M) = M.$ \\ \\
Theorem 2.2. For any module $ M_{R} $, the followings are equivalent:
\begin{enumerate}
\item[1.] M is rationally complete. 
\item[2.] For any right R-modules $ A \leq B $ such that $ Hom_{R}(B/A, E(M)) = 0 $, any R-homomorphism $ f : A \rightarrow M  $ can be extended to B.

\end{enumerate}

 If M is Rationally complete then M is not coretractable.  
We know that an R-module M is called monoform if every submodule is dense. Now we are ready to establish following proposition.\\ \\
Proposition 2.3. Let $ I \neq R $ be an ideal in a commutative ring R  and let $ M_{R}  = R/I $, then followings are equivqlent in M:
\begin{enumerate}
\item[1.] Every submodule in M is dense.
\item[2.] M is monoform.
\item[3.] M is not coretractable.
\item[4.] Every ideal in R/I is prime.
\end{enumerate}
Proof. $(1)\Rightarrow (2).$ \\
$ (1)\Rightarrow(3). $ \\
$ (1)\Leftrightarrow (4). $ Let $ M_{R} $ is a module and $ I\neq R $ is an ideal of R. Let us suppose that any submodule N of M is dense in M. Then $ Hom_{R}(N/I , M)=0 \Leftrightarrow f(I) \\ = 0 \Leftrightarrow f(ab) = 0 \Leftrightarrow ab =0 \Leftrightarrow either \ a= 0 \ or \ b = 0. 
\Leftrightarrow $ I is prime ideal.\\ \\
We say that M is a epi-retractable module if for every submodule N of M, we have $Hom(M, N) \neq 0. $ Our next result gives relation between epi-retractable module and coretractable module.\\ \\
Proposition 2.4.The following statements are equivalent for a module M:
\begin{enumerate}
\item[1.] M is epi-retractable.
\item[2.] There exist surjective homomorphisms $ M \rightarrow N \ and \  N \rightarrow M $ for some epi-retractable module N.
\item[3.]There exists a surjective homomorphism $ M/K \rightarrow M $ for some epi-retractable factor module M/K.
\item[4.] M is coretractable.
\end{enumerate}
Proof. $(1) \Rightarrow (2) \ and \ (2) \Rightarrow (3) $ By Lemma 2.7 \cite{10}. \\
$ (3) \Rightarrow (4) $ is obvious. \\
$ (4) \Rightarrow (1) $ Let N be any submodule of M. Since M is coretractable we have $ f:M/K\rightarrow M \neq0.$
So we have $ f_{\mid_{N/K}}:N/K \rightarrow N.$ Consider a projection map $ \pi: M/K \rightarrow  N/K $  and a canonical epimorphism $ g:M\rightarrow M/K.$ The composition $ f_{\mid_{N/K}} \pi g: M \rightarrow N $ implies that M is epi-retractable modules. \\ \\
Proposition 2.5. Let $ M_{R} $ be an epi-retractable module, then $ M_{R} $ is coretractable if and only if M is semiprimitive. \\ 
Proof. Let M is coretractable then M is semisimple (see {2.3}\cite{13}). So  M is semiprimitive . For the sufficient part let M is semiprimitive. Now since M is epi-retractable module, so we have a surjective mapping s.t. $ f(M,N) \\ \neq 0 ~\text{and}~ g(N,M) \neq 0 $ by {lemma 2.7} \cite{10} . Clearly we have a non zero homomorphism $ h(M/N, N) \neq 0. $ Thus the composition $ gh \neq 0.$ So M is coretractable. \\ \\
Proposition 2.6. Let M be a projective R-module, then every epi-retractable  module is coretractable.\\
Proof. Let M is a projective R-module and N is a submodule of M. Consider a non zero map $ f:M/N\rightarrow N. $  R is a right perfect ring so we have a surjective map $ \theta :P \rightarrow N ; $
P is projective cover of N.  
Since M is projective so $ f = \theta h \neq 0$ where $  h:M/N \rightarrow P $ is a non zero map. 
Since M is epi-retractable module so we have a surjective mapping from $  \pi: M\rightarrow N ~ \text{and} ~ g:N\rightarrow M. $ ( Prop. 2.1 \cite{10}).Clearly the mapping $ gf:M/K\rightarrow M \neq 0. $ Hence the proof.\\ \\
 An R-module M is satisfies (**) property if every $  f\epsilon End(M_{R}) $, $f$ is surjective mapping from M to M.\\
 We say that R is a right Kasch ring  if every simple right R- module M can be embedded in 
 $ R_{R}$. "Left Kasch ring" is defined similarly. R is called a Kasch ring if it is both right and left Kasch. A module M is called kasch module if it contains a copy of every simple module in $ \sigma (M). $\\ \\
Proposition 2.7.Let $ M_{R} $ be a epi-retractable, then
\begin{enumerate}

\item[1.] $ M_{R} $ is right kasch.
\item[2.] $ M_{R} $ satisfies (**) property.
\item[3.] $ M_{R} $ is coretractable.

\end{enumerate}
Proof. $ (1)\Rightarrow (2) $ Let $  f\epsilon End(M_{R}) ~ \text{and}~ f:M\rightarrow K $ be a non zero homomorphism where K is a proper submodule of M. Now M is right kasch so it contains copy of simple module. Thus we have $ f(M)=K=M. ~ f $ is  surjective. M satisfies (**) property.\\
$ (2)\Rightarrow (3)  $ Let M satisfies (**) property.Also Mis epi-retractable module so we have non-zero homomorphism from $ f:M/K \rightarrow N/K~\text{where}~K\leq N \leq M. $ Now we have a natural epimorphism $ f:N/K \rightarrow N $ and let $ g\epsilon End_{R}(M) \neq 0$ where $ g $ is an onto map from N to M. So we have a onto composition $ ghf:M/K \rightarrow M \neq 0. $  \\
$ (3)\Rightarrow (1) $  (Prop.2.14 \cite{13}). \\ \\
The following two proposition shows that coretractable modules satisfies (**) property and M is also a direct summand of a free epi-retractable R-modules.\\ \\
Proposition 2.8.  A Coretractable module  M satisfies (**) property.\\
Proof. Let M be a coretractable Module so we have a non zero homomorphism $ f:(M/K, M) \neq 0 $
for some proper submodule K. Consider a nonzero $ g\epsilon End(M_{R}) .$ We have $ gf:(M/K,M) \neq 0 $
and also the mapping $ gf $ is epimorphism. So M satisfies (**) property. \\ \\
Proposition 2.9. Let R be a ring and $ \alpha $ be an infinite ordinal $ \geq \mid R \mid.$ Suppose that
$ M = K \bigoplus N $ where K is a free epiretractable R-module with a basic set of cardinality $ \alpha $ and N is a $ \beta $-generated R-module with $ \alpha \geq \beta. \  Then, \ M_{R} $ is coretractable.\\ \\
Proof. K is epi-retractable so N is a homomorphic image of K. Since $ K \simeq K \bigoplus K,$ then there exist surjective homomorphisms $ f:M \rightarrow K $ and $g: K \rightarrow M.$ Let L = Kerf, then there exist natural homomorphism $ h:M/L \rightarrow M. $ Hence M is coretractable. \\ \\
A module M is said to be fully retractable if for every nonzero submodule N of M and every nonzero element $ g \epsilon Hom_{R}(N,M), $ we have  $ Hom_{R}(M,N)g  \neq 0.$ This definition motivated us to give our next result. \\ \\
Proposition 2.10. Every fully retractable module is coretractable if  $ Hom_{R}(M,N)g $ is epimorphism.\\ 
Proof. clear by prop. 2.4. \\ \\
Proposition 2.11.Let R be a ring. If $ M_{R} $ is cyclic CS module then M is coretractable. \\ 
Proof. since M is cyclic CS module then M is direct sum of uniform modules i.e.
M = $ \bigoplus^{n} _{i=1} N_{i}$ where each  $N_{i}'s $ is essential in M. Let K any submodule of M then K be essential in direct summand of M.
Now exploiting the fact that every quotient module of a cyclic module  is cyclic then we get for
$ Hom_{R}(M/K , M)\neq 0 $ giving us that M is coretractable.\\ \\
Proposition 2.12. Let R be a commutative artinian ring, then
\begin{enumerate}
\item[1.] If p be any prime ideal of R then p is maximal.
\item[2.] R is local ring.
\item[3.] $R_{R} $ is semisimple.
\item[4.] $ Hom_{R}(R/p, R) \neq 0. $
\item[5.] R is a kasch ring.
\item[6.] $ R_{R} $ is coretractable.
\end{enumerate}
Proof. Let p be any prime ideal of R, then R/p become an integral domain. i.e. R/p is a field consequently p is maximal ideal of R.\\
$ (1)\Rightarrow (2) $ Since any commutative ring R with prime ideal p has got unique maximal ideal.Thus R should be a local ring.\\
Since R is artinian, hence (3) is obvious.\\
$(1)\Rightarrow (4) $ Let p be a prime ideal of an Artinian ring R. Then R/p is a prime Artinian ring, but such rings are simple. Hence R is kasch.\\
$(5)\Rightarrow (6)$  Obvious. \\ \\
We know that direct sum of coretractable module is coretractable but in the case of direct product of coretractable module the result may not be same. Next result will validate our supposition.\\
Proposition 2.13. The infinite direct sum of coretractable module is not coretractable.\\ 
Proof.  Let $ K = K_{1}\oplus K_{2} \oplus ....\subset M $ be the infinite direct product of coretractable module where $ K_{i}, (i=1,2... )$ are projective. Consider a map \\ $ f:(M/K, K)~ \text{with}~ f(K)=0. $ 
Also suppose that $ M = M_{1}\oplus M_{2} . $ Since $ f(K)= 0 ~\text{so}~ f(M_{1})= 0 ~\text{and} ~ f(M_{2})= 0 \Rightarrow f(M)=0 \Rightarrow f:(M/K, M) = 0.$ Thus infinite direct sum of coretractable module is may not be coretractable. \\ \\
Proposition 2.14. Let R be hereditary ring, then $ M_{R} $ is coretractable.\\ 
Proof. Let N be a injective submodule of M. Then there exist a submodule K of M such that   
$ M = K \oplus N $ by  {lemma 2.8.} Let $ f(M/N, N) \neq 0. $ Consider a injection map $ g:N\rightarrow M. $ Clearly the composition 
$  gf (M/N,M) \neq 0.$ So $ M_{R} $ is coretractable.
\section{Mono coretractable Module}
Now we introduce Mono-coretractable module and we will discuss various propeties of mono-coretractable modules. \\
 An R-Module M is said to be Mono-coretractable if  $ Hom_{R}(M/K,M)\neq 0 $ is monomorphism where K is a submodule of M. A coretractable  modules are not mono-coretractable in general for example it can be easily seen that $ Z_{4} ~\text{over} ~ Z $ is coretractable but not mono-coretractable.\\ \\
Let us discuss in this section various properties of mono-coretractable modules.\\
Proposition 3.1. For a module $ M_{R} $ following is equivqlent: 
\begin{enumerate}
\item[1.] M is mono coretractable.
\item[2.] Every non zero partial endomorphism of M/K  is monomorphism.
\end{enumerate}
Proof. $(1)\Rightarrow (2) $ Let M is mono coretractable.We have $ HOM_{R}(M/K ,M)\neq 0. $ 
Take $ N/K \leq M/K \ where \ N \leq K \leq M $ and let $  f:N/K\rightarrow M/K $ be a nonzero partial endomorphism of M/K. Then we have  $ f:N/K\rightarrow N.kerf\cong imf. \ Also \ f_{1}(N/kerf) \rightarrow M/K  $ is monomorphism.\\
 The composition $ N/K \rightarrow M/K \rightarrow M  $ is monomorphism. Thus kerf=K.  Thus Every non zero partial endomorphism of M/K  is monomorphism.\\
 $(1)\Rightarrow (2) $  let Every non zero partial End(M/K)  is monomorphism and let M be module s.t.
 $  f:M/K \rightarrow M \neq 0 . \ \text{Now} \ f:N/K\rightarrow M/K   $ is monomorphism. Consider a  projection map $ h:N/K \rightarrow M $ define as $ f(x +K) = x , x \epsilon N.  $ Clearly h is monomorphism.\\
 Since h is projection map, so h = fg and the composition $  N/K \rightarrow M/K \rightarrow M $ should be monomorphism  $ \Rightarrow $ f is monomorphism $ \Rightarrow $  M is mono coretractable. \\ \\
 Proposition 3.2. Direct sum of mono coretractable modules is mono coretractable.\\ \\
Proof. Obvious. \\ \\
Proposition 3.3. Let $ M_{R} $  be a non singular uniform module in R then $ M_{R} $ is not mono coretractable. \\ \\
Proof. take any $ N\leq M $ then M/N is singular. Hence $ Hom_{R}(M/N , M)=0.$
So M is not mono coretractable. \\ \\
Now we discuss some condition for mono-coretractability.\\
Proposition 3.4. For a ring R, Following holds:
\begin{enumerate}
 \item[1.] R is a right Kasch ring.
\item[2.] $ R_{R} $ is a mono- coretractable module.
\item[3.] Every finitely generated free right R-module is mono-coretractable.
\item[4.] $R_{R}$ has no proper dense submodules.
\item[5.] $R_{R}$  not monoform.
\end{enumerate}
Proof:-$(1)\Rightarrow(2) $ For any maximal ideal I in R, $Hom_{R}(R/I, R) \neq 0.$
 Thus $R_{R} $ is monocoretractable. \\
 $ (2)\Rightarrow (3). $ Obvious\\
$(2) \Rightarrow(4) $ Take a R-module M and Consider $ K\leq N \leq M. $ Due to Mono coretractability  $Hom_{R}(M/K, M) \neq 0\Rightarrow Hom_{R}M/K, N) \neq 0.$ So M has no propoer dense submodules. \\
$(4) \Rightarrow(5). $  Obvious.



\begin{thebibliography}{}
\bibitem{1} Gooderl, K.R. and Warfield R.B. Jr, \emph{An Introduction to non-commutative Noetherian Rings}, second edition, Vol. 61, London Mathematical Society, Student Texts,Cambridge Press,Cambridge, (2004).
\bibitem{2} Mc connel, J.C. and Robson, \emph{non-commutative Noetherian Rings}, Graduate studies in Mathematics, Vol. 30, American Mathematical Society, (1987).
\bibitem{3}S.M. Khuri,Endomorphism rings and lattice isomorphism,\emph{J. Algebra,} {\bf 56}(2),401-408,(1979).
\bibitem{4} S.M. Khuri,Endomorphism rings of nonsingular modules, \emph{Ann. Sci. Math. Quebec,}{\bf 4}(2),145-152,(1980).
\bibitem{5}P.F. Smith. Modules With Many Homomorphism,\emph{Journal of pure and Applied Algebra}, {\bf 197}, 305-321,(2005).
\bibitem{6}P.F. Smith and M.R. Vedadi, Essentially Compressible Modules and Rings, \emph{Journal of Algebra}, {\bf 304}, 812-831,(2006).
\bibitem{7}Abhay K. Singh, Essentially Slightly Compressible Modules and Rings,\emph{Asian-European Journal of Mathematics(World Scientific)}, {\bf 5}(2),(2012).
\bibitem{8}Virginia Silva Rodrigues and Alveri Alves Sant’Ana, A note on a Problem due to Zelmanowitz, \emph{Algebra and Discrete Mathematics}, No.3, 85-93,(2009).
\bibitem{9}M.R. Vedadi, Essentially Retractable Modules, Journal of Science, \emph{Islamic Republic of Iran} {\bf 18}(4), 355-360,(2007).
\bibitem{10}A. Ghorbani and M.R. Vedadi, Epi-retractable Modules and some Applications, \emph{Bulletin of Iranian mathematical Society,} {\bf 35} (1), 155-166, (2009).
\bibitem{11}J. Zelmanowitz, An extension of the Jacobson density theorem, \emph{Bulletin American Mathematical Society},{\bf 82}(4),551-553, (1976).
\bibitem{12}J. Zelmanowitz, Weakly primitive rings,\emph{Communication of Algebra}, {\bf 9}(1),(23-45), (1981).

\bibitem{13} B. Amini, M. Ershad and H. Sharif, Coretractable Modules, \emph{Journal of  Australian Mathematical Society} {\bf 86}, 289-304, (2009).

\bibitem{14} J. Zemlicka, Completely Coretractable Rings, \emph{Bulletin of the Iranian Mathematical Society} {\bf 39}(3),523-528,(2013).
\bibitem{15} A.M.A. Al-Hosainy and H.A. Kadhim,Co-compressible Module,\emph{International Research Journal of Scientific Findings}{\bf 1}(6),(2014).

\bibitem{16}S.M. Khuri,Endomorphism rings of a nonsingular retractable modules,\emph{East-West J. Math.}
{\bf 2}(2).



%
\end{thebibliography}


\end{document}